\documentclass[a4paper,11pt]{article}

  \renewcommand\appendix{\par
  \setcounter{section}{0}
  \setcounter{subsection}{0}
  \setcounter{figure}{0}
  \setcounter{table}{0}
  \renewcommand\thesection{ Appendix \Alph{section}}
  \renewcommand\thefigure{\Alph{section}\arabic{figure}}
  \renewcommand\thetable{\Alph{section}\arabic{table}}
}

 \usepackage{amsmath, amsthm, amssymb}
\usepackage{amssymb}
\usepackage{graphicx}
\usepackage{algorithmic}
\usepackage{algorithm}

\usepackage{tikz}
\usetikzlibrary{calc}
\usetikzlibrary{patterns}
\tikzstyle{mybox} = [draw=black, fill=white,  thick,
    rectangle, inner sep=10pt, inner ysep=20pt]
\tikzstyle{mybox} = [draw=black, fill=white,  thick,
    rectangle, inner sep=2pt, inner ysep=2pt]

\usepackage[margin=1.2in]{geometry}

\usepackage{setspace}

\newtheorem{definition}{Definition}[section]

\newtheorem{theorem}{Theorem}[section]

\newtheorem{remark}{Remark}[section]
\newtheorem{example}{Example}[section]

\begin{document}

\begin{center}
	{\Large \bf Generalized Kneser hypergraphs and a new variant of colorability defect}\\
	\vspace{0.3 cm}
	{\bf Sharareh Alipour$^\dag$, Amir Jafari$^\ddag$, and Saeed Shaebani$^\S$\\
		{\it $^\dag$ Department of Software Engineering}\\
		{\it New Uzbekistan University}\\
		{\it Tashkent, Uzbekistan}\\
		{\tt sh.alipour@newuu.uz}\\
		{\it $^\ddag$ Department of Mathematics}\\
		{\it New Uzbekistan University}\\
		{\it Tashkent, Uzbekistan}\\
		{\tt a.jafari@newuu.uz}\\
		{\it $^\S$ School of Mathematics and Computer Science}\\
		{\it Damghan University}\\
		{\it Damghan, Iran}\\
		{\tt shaebani@du.ac.ir}\\
		
	}
\end{center}

\begin{abstract}
\noindent
Dol'nikov and K{\v{r}}{\'{\i}}{\v{z}} introduced the notion of $r$-colorability defect of hypergraphs in order to provide a strong and sharp lower bound for the chromatic number of general Kneser hypergraphs ${\rm KG}^r({\mathcal F})$.

\noindent Given a hypergraph ${\mathcal F}$ together with a nonnegative integer $s$ such that $s<|F|$ for all hyperedges $F$ in $E({\mathcal F})$, one defines a more general Kneser hypergraph ${\rm KG}^r({\mathcal F}, s)$ with vertex set $E({\mathcal F})$ whose hyperedge set consists of all $r$-subsets $\{A_1, A_2, \dots , A_r\}$ such that $|A_i\cap A_j|\le s$ for all $1\le i<j\le r$.

\noindent In this paper, we generalize Dol'nikov and K{\v{r}}{\'{\i}}{\v{z}}'s result by introducing the notion of complementary colorability defect of hypergraphs and providing a lower bound for $\chi ({\rm KG}^p({\mathcal F} , s))$ for prime numbers $p$ by means of the celebrated ${\mathbb Z_p}$-Tucker Lemma. Our result, also extends Daneshpajouh's lower bound for $\chi ({\rm KG}^r(n , k , s))$.\\

\end{abstract}

\section{Introduction}

A {\it hypergraph} $\mathcal F$, is an ordered pair $(V,E)$, usually shown by
$(V({\mathcal F}),E({\mathcal F}))$ instead, where $V$ is a finite set called the {\it vertex set} of $\mathcal F$, and $E$ which is called the {\it hyperedge set} of ${\mathcal F}$, is a subset of ${\mathcal P} (V({\mathcal F})) \setminus \{\varnothing\}$. Indeed, a hypergraph consists of a finite vertex set $V$ equipped with a family of some nonempty subsets of $V$.
By an {\it $r$-uniform} hypergraph, we mean a hypergraph with $|F|=r$ for each hyperedge $F$ in $E({\mathcal F})$. In particular, any 2-uniform hypergraph is called simply a
{\it graph}. Throughout, for any nonnegative integer $n$, we abbreviate $\{1,2,3, \dots , n\}$ by $[n]$. Also, for any set $A$ and any nonnegative integer $k$, we define ${A \choose k}$ equal to the set of all $k$-subsets of $A$, that is,
${A \choose k} := \{B : B\subseteq A \ {\rm and} \ |B|=k\}.$ 

Let ${\mathcal F}$ be a hypergraph. By a {\it vertex coloring} of ${\mathcal F}$, we mean a function
$c: V({\mathcal F}) \longrightarrow C$ such that $|c(F)| > 1$ for each hyperedge $F$ of ${\mathcal F}$. The {\it chromatic number} of ${\mathcal F}$, denoted by $\chi ({\mathcal F})$, is defined as the minimum size of a set $C$ that a vertex coloring $c: V({\mathcal F}) \longrightarrow C$ exists. If ${\mathcal F}$ has a hyperedge $F$ with $|F|=1$, then we define $\chi ({\mathcal F}) := \infty$.

\noindent Let ${\mathcal F}$ be a hypergraph and $r$ be an integer greater than or equal to 2.
The symbol ${\rm KG} ^{r} ({\mathcal F})$ stands for an $r$-uniform hypergraph with vertex set
$E({\mathcal F})$; and an $r$-subset $\{F_{1} , F_{2} , \dots , F_{r}\}$ of $V \left( {\rm KG} ^{r} ({\mathcal F}) \right)$ belongs to $E \left( {\rm KG} ^{r} ({\mathcal F}) \right)$ iff
$F_{i}  \cap F_{j} = \varnothing$ for all distinct $i$ and $j$ in $[r]$.

\noindent Let $n, k,$ and $r$ be positive integers with $r\geq 2$, and consider the {\it complete $k$-uniform hypergraph} $K_{n}^{k} := \Big( [n], {[n] \choose k} \Big)$. The graph
${\rm KG} ^{2} \left( K_{n}^{k} \right)$ is called the {\it Kneser graph} ${\rm KG} (n,k)$, and the hypergraph ${\rm KG} ^{r} \left( K_{n}^{k} \right)$ is known as the {\it usual Kneser hypergraph} ${\rm KG}^{r} (n,k)$.
Kneser \cite{Kneser1978} conjectured that $\chi \left( {\rm KG} ^{2} \left( K_{n}^{k} \right) \right) = n-2k+2$, provided that $n\geq 2k-1$. This conjecture was settled by Lov\'{a}sz \cite{Lovasz1978}. Also, Alon, Frankl, and Lov\'{a}sz \cite{AlonFranklLovasz1986}, by solving a more general conjecture of Erd\"{o}s \cite{Erdos1973}, extended
this result by showing that whenever $n \geq r(k-1)+1$, then
$$\chi \left( {\rm KG} ^{r} \left( K_{n}^{k} \right) \right) = \left\lceil  \frac{n-r(k-1)}{r-1}  \right\rceil .$$

For a hypergraph $\mathcal F$, Dol'nikov \cite{Dol'nikov1988} defined the notion of {\it 2-colorability defect} ${\rm cd}^2(\mathcal F)$, as the smallest size of a subset $X_0$
of $V({\mathcal F})$ so that there is a partition $\{X_1 , X_2\}$ of $V({\mathcal F}) \setminus X_0$ where there is no $F$ in $E({\mathcal F})$ with $F\subseteq X_i$ for $i=1$ or $i=2$. For example, if $n \geq 2k-1$, then ${\rm cd}^2(K_{n}^k)=n-2k+2$. Dol'nikov \cite{Dol'nikov1988} proved that the chromatic number of ${\rm KG}^2(\mathcal F)$ is at least ${\rm cd}^2(\mathcal F)$ for each hypergraph ${\mathcal F}$ in general.

\begin{theorem} \cite{Dol'nikov1988} \label{Dol'nikov}
	Every hypergraph ${\mathcal F}$ satisfies
	$\chi({\rm KG}^2({\mathcal F} ))\ge  {\rm cd}^2({\mathcal F} ).$
\end{theorem}

The notion of 2-colorability defect was extended by K{\v{r}}{\'{\i}}{\v{z}} to {\it $r$-colorability defect} ${\rm cd}^r(\mathcal F)$ for any integer $r\ge 2$ to be the smallest size of a subset $X_0$ of $V({\mathcal F})$ so that there is a partition $V({\mathcal F}) \setminus X_0 = X_1\cup\dots\cup X_r$ where there is no $F$ in $E({\mathcal F})$ with $F\subseteq X_i$ for some $i$ in $\{1,2,\dots , r\}$. For example, if $n \geq r(k-1)+1$ then ${\rm cd}^r(K_{n}^k)=n-r(k-1)$. The nontrivial part of the Erd\"os' conjecture was that the chromatic number of ${\rm KG}^r(K_{n}^k)$ is at least $\left\lceil \frac{{\rm cd}^r(K_{n}^k)}{r-1}\right\rceil$ which was proved by Alon, Frankl, and Lov\'asz \cite{AlonFranklLovasz1986}. K{\v{r}}{\'{\i}}{\v{z}} proved in \cite{Kriz, Kriz2000} that more generally, the chromatic number of ${\rm KG}^r(\mathcal F)$ is at least $\left\lceil \frac{{\rm cd}^r(\mathcal F)}{r-1}\right\rceil$.

\begin{theorem} \cite{Kriz, Kriz2000} \label{Kriz}
	If ${\mathcal F}$ is a hypergraph and $r$ is an integer with $r\geq 2$, then we have
	$$\chi({\rm KG} ^{r} ({\mathcal F}))\geq \left\lceil \frac{{\rm cd}^r({\mathcal F})}{r-1}\right\rceil.$$
\end{theorem}

For positive integers $n$ and $k$, and a nonnegative integer $s$ that satisty $n \geq k > s$, Frankl \cite{Frankl1985} defined the graph ${\rm KG}^2(K_{n}^k,s)$ to be a graph whose vertex set is ${[n] \choose k}$, and $\{F_1 , F_2\}$ forms an edge whenever $|F_1\cap F_2|\le s$. This definition was extended by Daneshpajouh \cite{Dan} to an $r$-uniform hypergraph ${\rm KG} ^ r (K_{n}^k,s)$ whose vertex set is ${[n] \choose k}$, and $\{F_1 , F_2 , \dots , F_r\}$ forms a hyperedge whenever $|F_i \cap F_j| \leq s$ for any two distinct $i$ and $j$ in $[r]$.

\noindent
It was shown by Daneshpajouh in \cite{Dan} that the chromatic number of ${\rm KG}^p(K_{n}^k, s)$ is at least $\left\lceil \frac{n-p(k-s-1)}{p-1}\right\rceil$
for all prime numbers $p$. 

\begin{theorem} \label{DanTh} \cite{Dan}
	Let $p$ be a prime number.
	If $n, k,$ and $s$ are integers such that $k > s \geq 0$ and $n\geq p(k-1)+1$, then we have
	$$\chi \left( {\rm KG} ^p \left( K_{n}^k,s \right) \right) \geq \left\lceil \frac{n-p(k-s-1)}{p-1} \right\rceil .$$ 
\end{theorem}

For a hypergraph ${\mathcal F}$ and any integers $r$ and $s$ which satisfy $r\ge 2$ and also $0 \leq s < |F|$ for all hyperedges $F$ of ${\mathcal F}$, one defines the {\it generalized Kneser hypergraph} ${\rm KG}^r({\mathcal F}, s)$ to be an $r$-uniform hypergraph with vertex set $V \left(  {\rm KG}^r({\mathcal F}, s)  \right) := E({\mathcal F})$; and $\{F_1, F_2 , \dots , F_r\}$ represents a hyperedge of ${\rm KG}^r({\mathcal F}, s)$ iff $|F_i\cap F_j|\le s$ for all integers $i$ and $j$ with $1 \le i < j \le r$.

In order to obtain a lower bound for ${\rm KG}^r(\mathcal F, s)$,
Azarpendar and Jafari \cite{Azarpendar} made the following definition.
\begin{definition}
	For a hypergraph $\mathcal F$ and any integers $r$ and $s$ which satisfy $r\ge 2$ and also $0 \leq s < |F|$ for all hyperedges $F$ of $\mathcal F$, the {\it general $r$-colorability defect} ${\rm cd}^r(\mathcal F,s)$ is the smallest size of a subset $X_0$ of $V({\mathcal F})$ so that there is a partition $V({\mathcal F}) \setminus X_0 = X_1\cup \dots\cup X_r$ where there is no $F$ in
	$E({\mathcal F})$ such that $|F \setminus X_i|\le s$ for some $i$ in $\{1,2, \dots, r\}$.
\end{definition}
It is clear that ${\rm cd}^r(\mathcal F,0)={\rm cd}^r(\mathcal F)$, since $|A \setminus X_i|\le 0$ if and only if $A\subseteq X_i$.
As an improvement of Dol'nikov and K{\v{r}}{\'{\i}}{\v{z}}'s Theorems, it was proved in \cite{Azarpendar} that the chromatic number of ${\rm KG}^r(\mathcal F)$ is at least $\left\lceil \frac{{\rm cd} ^r \left({\mathcal F} , \left\lfloor \frac{s}{2} \right\rfloor \right) }{r-1} \right\rceil$ for any positive integer $r$ greater than or equal to 2.
\begin{theorem} \label{ecd} \cite{Azarpendar}
	Let $\mathcal F$ be a hypergraph and $r$ and $s$ be some integers that satisfy $r\ge 2$ and also $0 \leq s < |F|$ for all hyperedges $F$ of $\mathcal F$. Then we have
	$$\chi \left( {\rm KG}^r ({\mathcal F} , s)  \right) \geq
	\left\lceil \frac{{\rm cd} ^r \left({\mathcal F} , \left\lfloor \frac{s}{2} \right\rfloor \right) }{r-1} \right\rceil.$$
\end{theorem}

\noindent
Azarpendar and Jafari \cite{Azarpendar} conjectured that Theorem \ref{ecd} remains true if one replaces $\left\lfloor \frac{s}{2} \right\rfloor$ with $s$. Shaebani \cite{Shaebani2022}
gave a negative answer to this conjecture by showing that $\left\lfloor \frac{s}{2} \right\rfloor$ is the best possible and no integer greater than $\left\lfloor \frac{s}{2} \right\rfloor$ could be replaced.

\noindent
We know that if $n \geq r(k-1) + 1$, then 
$${\rm cd}^r \left(K_{n}^k, \left\lfloor \frac{s}{2} \right\rfloor \right)=n-r\left(k-\left\lfloor \frac{s}{2} \right\rfloor-1\right);$$
and therefore, the assertion of Theorem \ref{ecd} turns into
$$\chi \left( {\rm KG}^r \left(K_{n}^k , s \right)  \right) \geq
\left\lceil \frac{n-r\left(k-\left\lfloor \frac{s}{2} \right\rfloor-1\right) }{r-1} \right\rceil.$$
So, the lower bound in Theorem \ref{ecd} is not quite as strong as the result in Theorem \ref{DanTh}.

The goal in this article is to define a new variant of colorability defect that could be used to provide a lower bound for $\chi \left( {\rm KG}^r ({\mathcal F} , s)  \right)$ in such a way that Theorem \ref{DanTh} would be derived as a straightforward consequence. We call this new invariant {\it the $s$-th complementary $r$-colorability defect} of $\mathcal F$ and denote it by ${\rm ccd}^r(\mathcal F,s)$; as follows.
\begin{definition}
	Let $\mathcal F$ be a hypergraph, $r$ be an integer with $r\geq 2$, and
	$s$ be a nonnegative integer for which $s < |F|$ for each hyperedge $F$ in $E({\mathcal F})$.
	The {\it $s$-th complementary $r$-colorability defect} of ${\mathcal F}$, denoted by ${\rm ccd}^r({\mathcal F} , s)$, is the smallest size of a subset $X_0$ of $V({\mathcal F})$ such that there is a partition $V({\mathcal F}) \setminus X_0 = X_1\cup \dots\cup X_r$ where there is no $F$ in
	$E({\mathcal F})$ with $|F\cap X_0|\le s$ and $F\subseteq X_0\cup X_i$ for some $i \in \{1, 2, \dots , r\}$.
	Roughly speaking, ${\rm ccd}^r({\mathcal F} , s)$ is defined as the minimum size of a subset $X_0$ of $V( \mathcal F)$ that a partition $\{X_1, X_2, \dots , X_r\}$ of $V({\mathcal F}) \setminus X_0$ exists in such a way that for each hyperedge $F$ in $E ({\mathcal F})$ and each $i$ in $[r]$, at least one of the following two conditions does not hold:
	\begin{itemize}
		\item $F \subseteq X_0 \cup X_i ;$
		\item $|F \cap X_0| \leq s.$
	\end{itemize}
	In this definition, it is necessary to clarify that in a partition $V({\mathcal F}) \setminus X_0 = X_1 \cup X_2 \cup \dots \cup X_r$, it is permitted that $X_{i}=\varnothing$ for some $i$ in $\{0,1,2,\dots,r\}$, and there is no constraint for $X_{i}$'s to be nonempty.
\end{definition}

\noindent For $s=0$, the $s$-th complementary $r$-colorability defect equals the classical $r$-colorability defect, ${\rm cd}^r ({\mathcal F})$, which is the minimum number of vertices of ${\mathcal F}$ that removing them from ${\mathcal F}$ yields a remaining $r$-colorable subhypergraph of ${\mathcal F}$.

\noindent
We always have ${\rm ccd}^r(\mathcal F,s)\le {\rm cd }^r(\mathcal F,s)$. This inequality may become strict with a wide gap but there are also examples where it becomes an equality.
Here are two relevant examples.
\begin{example}
	Let $s:=1$ and $r \geq 2$ and ${\mathcal F}$ be the {\it star graph} $K_{1,n}$ with vertex set $V({\mathcal F}) := [n+1]$ and edge set
	$E({\mathcal F}) :=\{\{i,n+1\} : \ i \in [n]\}$.
	Now, by defining $X_{0} = \dots = X_{r-2} = \varnothing$ and
	$X_{r-1} = [n]$ and $X_{r} = \{n+1\}$, we observe that ${\rm ccd}^r(\mathcal F,s) = 0$ while ${\rm cd }^r(\mathcal F,s) = n+1$.
\end{example}

\begin{example}
	Let $\{n,r,k\} \subseteq \{2,3,4, \dots\}$. Now,
	 if $n \geq r(k-1)+1$ and $0\le s < k$, then 
	$${\rm cd}^r \left(K_{n}^k, s \right)={\rm ccd}^r \left(K_{n}^k, s \right)=n-r\left(k-s-1\right).$$
\end{example}

The following theorem is the main result of this paper.
\begin{theorem} \label{Saeed}
	Let $\mathcal F$ be a hypergraph, $p$ be a prime number, and $s$ be a nonnegative integer for which $s < |F|$ for all hyperedges $F$ of ${\mathcal F}$. Then, we have
	$$\chi({\rm KG}^p({\mathcal F} , s))\ge \left\lceil \frac{{\rm ccd}^p({\mathcal F} , s )}{p-1}\right\rceil.$$
\end{theorem}
\noindent Since the assertion of Dol'nikov and K{\v{r}}{\'{\i}}{\v{z}}'s Theorems
could be stated as 
$$\chi({\rm KG}^r({\mathcal F} , 0))\ge \left\lceil \frac{{\rm ccd}^r({\mathcal F} , 0 )}{r-1}\right\rceil,$$
Theorem \ref{Saeed} is a natural extension of them for all prime numbers $r$.

\begin{remark}
	If $p$ is a prime number and integers $n, k,$ and $s$ satisfy $k > s \geq 0$ and $n\geq p(k-1)+1$, then by considering
	the hypergraph ${\mathcal F} := K_{n}^k$ we would have ${\rm ccd}^p({\mathcal F} , s )= n-p(k-s-1)$. So, as a corollary of Theorem \ref{Saeed}, the Theorem \ref{DanTh} of Daneshpajouh follows immediately.
\end{remark}

\section{Preliminaries and The Proof}

In this section, we will be concerned with proving Theorem \ref{Saeed}.
The main idea of the proof is to apply a celebrated powerful theorem of Ziegler and Meunier, which is called The ${\mathbb Z_p}$-Tucker Lemma \cite{Meunier2011,Ziegler2002,Ziegler2006}.

\noindent Let $p$ be a prime number and ${\mathbb Z_p}:=\{\omega , \omega ^2 , \omega ^3 , \dots , \omega ^ p\}$ be the multiplicative cyclic group of order $p$ with a generator $\omega$, and with the additional property $0 \notin {\mathbb Z_p}$. Every $n$-tuple $X:=(x_1 , x_2 , \dots , x_n)$ in the $n$ times Cartesian product $({\mathbb Z_p} \cup \{0\})^n$ could be uniquely correspond to a $p$-tuple $(X_1 , X_2 , \dots , X_p)$, where $X_i:= \{j : \ j \in [n] \ {\rm and}\ x_j=\omega ^i\}$ for each $i$ in $[p]$. Also, we define $X_0:=\{j : \ j \in [n] \ {\rm and} \ x_j=0\}$.

\noindent If we agree that for each $\omega ^i$ in ${\mathbb Z_p}$, the relation $$\omega ^i . 0 = 0.0 = 0$$ holds, then for each $X:=(x_1 , x_2 , \dots , x_n)$ in $({\mathbb Z_p} \cup \{0\})^n$
we can define $$\omega ^i . X := (\omega ^i . x_1 \ , \ \omega ^i . x_2  \ , \  \dots  \ , \  \omega ^i . x_n).$$
\begin{theorem} {\rm (The Ziegler-Meunier ${\mathbb Z_p}$-Tucker Lemma  \cite{Meunier2011,Ziegler2002,Ziegler2006})}
	
	\noindent Let $m$ and $n$ be positive integers, $\alpha$ be a nonnegative integer, $p$
	be a prime number, and $m \geq \alpha$. Also, let
	${\mathbb Z_p}=\{\omega , \omega ^2 , \omega ^3 , \dots , \omega ^ p\}$ be the multiplicative cyclic group of order $p$ with a generator $\omega$ for which $0 \notin {\mathbb Z_p}$, and
	$$\lambda: ({\mathbb Z_p} \cup \{0\}) ^ n \setminus \{(0,0, \dots , 0)\} \longrightarrow \mathbb Z_p\times [m]$$
	$$\lambda(X)=(\lambda_1(X),\lambda_2(X))$$
	be a map satisfying the following three conditions:
	\begin{itemize}
		\item For each $X$ in $({\mathbb Z_p} \cup \{0\}) ^ n \setminus \{(0,0, \dots , 0)\}$
		and each $\omega ^i$ in ${\mathbb Z_p}$ we have
		$\lambda (\omega ^i.X)=(\omega ^i.\lambda_1(X),\lambda_2(X))$. In this case, $\lambda$
		is called a ${\mathbb Z_p}$-equivariant map.
		\item If $X$ and $Y$ are arbitrary elements in $({\mathbb Z_p} \cup \{0\}) ^ n \setminus \{(0,0, \dots , 0)\}$ in such a way that $\lambda_2(X) = \lambda_2(Y) \leq \alpha$ and $X_{i} \subseteq Y_{i}$ for each $i$ in $[p]$, then we have
		$\lambda_1(X) = \lambda_1(Y)$.
		\item If $X^{(1)} , X^{(2)} , \dots , X^{(p)}$ are arbitrary
		elements in $({\mathbb Z_p} \cup \{0\}) ^ n \setminus \{(0,0, \dots , 0)\}$
		in such a way that $\lambda_2 \left(X^{(1)}\right) = \lambda_2 \left(X^{(2)}\right) = \dots = \lambda_2 \left(X^{(p)}\right) \geq \alpha + 1 $ and $X^{(1)}_{i} \subseteq X^{(2)}_{i} \subseteq X^{(3)}_{i} \subseteq \dots \subseteq X^{(p)}_{i} $ for each $i$ in $[p]$, then  $\lambda_1 \left(X^{(1)}\right) , \lambda_1 \left(X^{(2)}\right) , \dots , \lambda_1 \left(X^{(p)}\right)$ are not pairwise distinct, that is,
		$\left| \left\{ \lambda_1 \left(X^{(1)}\right) , \lambda_1 \left(X^{(2)}\right) , \dots , \lambda_1 \left(X^{(p)}\right) \right\} \right| \leq p-1.$
	\end{itemize}
	Then, we have
	$$\alpha + (m-\alpha)(p-1) \geq n.$$	
\end{theorem}

\subsection{The Proof of Theorem \ref{Saeed}}

Now, we are in a position to prove Theorem \ref{Saeed}, as follows.

\begin{proof}
	Without loss of generality, we may assume that $V({\mathcal F})=[n]$.
	Let us regard an arbitrary total ordering $\preccurlyeq$ on ${\mathcal P} [n]$, the power set of $[n]$, such that $A\subseteq B$ implies $A \preccurlyeq B$. Also, let
	$t:=\chi({\rm KG}^p({\mathcal F} , s))$ and the function $c: E({\mathcal F}) \rightarrow [t]$
	be a proper vertex coloring of ${\rm KG}^p({\mathcal F} , s)$. Put $\alpha := n-{\rm ccd}^p({\mathcal F} , s )$ and $m:=\alpha + t$.
	Also, let the map
	$$\lambda: ({\mathbb Z_p} \cup \{0\}) ^ n \setminus \{(0,0, \dots , 0)\} \longrightarrow \mathbb Z_p\times [m]$$
	$$\lambda(X)=(\lambda_1(X),\lambda_2(X))$$
	be defined as follows:
	\begin{itemize}
		\item If $X=(X_1, X_2, \dots , X_p) \in ({\mathbb Z_p} \cup \{0\}) ^ n \setminus \{(0,0, \dots , 0)\}$ and $\displaystyle\sum_{j=1}^p \left| X_j \right| \geq \alpha + 1$,
		then $\left| X_0 \right| = n - \displaystyle\sum _{j=1} ^p \left| X_j \right| \leq {\rm ccd}^p({\mathcal F} , s ) - 1$. Hence
		$\left| X_0 \right| < {\rm ccd}^p({\mathcal F} , s )$; and therefore there exists a hyperedge 
		$F$ in $E({\mathcal F})$
		and a number $i$ in $[p]$ for which both conditions $F \subseteq X_0 \cup X_i$ and $\left| F \cap X_0 \right| \leq s$ hold.
		Let $\tilde{F}$ be the maximum of all such hyperedges $F$ with respect to the total ordering $\preccurlyeq $, accompanied with some element $\tilde{i}$ in $[p]$ for which
		$\tilde{F} \subseteq X_0 \cup X_{\tilde{i}}$ and $\left| \tilde{F} \cap X_0 \right| \leq s$.
		On one hand, since $\tilde{F} \subseteq X_0 \cup X_{\tilde{i}}$ and
		$\left| \tilde{F} \setminus X_0 \right| = \left| \tilde{F} \right| -  \left| \tilde{F} \cap X_0 \right| \geq \left| \tilde{F} \right| - s >0$,
		we have $\tilde{F}\cap X_{\tilde{i}} \neq \varnothing$.
		On the other hand, because of $\tilde{F} \subseteq X_0 \cup X_{\tilde{i}}$ we must have $\tilde{F}\cap X_j = \varnothing$ for all other numbers $j$ in $[p] \setminus \Bigl\{\tilde{i}\Bigr\}$.
		Accordingly, such a number ${\tilde{i}}$ would absolutely depend on $\tilde{F}$ and it is certainly unique with respect to $\tilde{F}$. So, it could be easily called ${\rm index}\left(\tilde{F}\right)$. Now, for determining $\lambda (X)$, put $$\lambda (X) := \left(\omega ^{{\rm index}\left(\tilde{F}\right)} \ , \ \alpha + c \left( \tilde{F} \right) \right).$$
		\item If $X=(X_1, X_2, \dots , X_p) \in ({\mathbb Z_p} \cup \{0\}) ^ n \setminus \{(0,0, \dots , 0)\}$ and $\displaystyle\sum_{j=1}^p \left| X_j \right| \leq \alpha $, then let
		$X_i := \max \{X_1 , X_2 , \dots , X_p\}$ with respect to the total ordering $\preccurlyeq$. Now, put
		$$\lambda (X) := \left(\omega ^{i} \ , \ \sum_{j=1}^p \left| X_j \right|  \right).$$
	\end{itemize}

The map $\lambda$ is ${\mathbb Z_p}$-equivariant. Also, let $X$ and $Y$ be two arbitrary elements in
$({\mathbb Z_p} \cup \{0\}) ^ n \setminus \{(0,0, \dots , 0)\}$ for which $\lambda _2 (X) = \lambda _2 (Y) \leq \alpha$ and $X_j \subseteq Y_j$ for each $j$ in $[p]$. Now, on one hand, according to the definition of $\lambda$, we have:
$$\begin{array}{llcrr}
	\lambda_2 (X) = \displaystyle\sum_{j=1}^p \left| X_j \right| & & {\rm and} & & \lambda_2 (Y) = \displaystyle\sum_{j=1}^p \left| Y_j \right| .
\end{array}$$
On the other hand, since $X_j \subseteq Y_j$ for each $j$ in $[p]$,
we have:
$$\left| X_j \right| \leq \left| Y_j \right| \mbox{ for all } j \mbox{ in } [p] .$$
Thus, because of $\lambda_2 (X) = \lambda_2 (Y)$, one obtains that $X_j = Y_j$ far each $j$ in $[p]$. Therefore $X=Y$ and in particular $\lambda _1 (X) = \lambda _1 (Y)$.

Finally, in order to show that $\lambda$ satisfies all the hypotheses of the celebrated ${\mathbb Z_p}$-Tucker Lemma, let $X^{(1)} , X^{(2)} , \dots , X^{(p)}$ be arbitrary
elements in $({\mathbb Z_p} \cup \{0\}) ^ n \setminus \{(0,0, \dots , 0)\}$
such that $X^{(1)}_{i} \subseteq X^{(2)}_{i} \subseteq X^{(3)}_{i} \subseteq \dots \subseteq X^{(p)}_{i} $ for each $i$ in $[p]$, and also,
$$\lambda_2 \left(X^{(1)}\right) = \lambda_2 \left(X^{(2)}\right) = \dots = \lambda_2 \left(X^{(p)}\right) \geq \alpha + 1 .$$
Therefore, on account of the definition of $\lambda$, there exist some hyperedges
$\tilde{F_1} , \tilde{F_2} , \dots , \tilde{F_p}$
in $E({\mathcal F})$ such that
$$\lambda_2 \left(X^{(1)}\right) = \alpha + c\left(\tilde{F_1}\right) \ , \ \lambda_2 \left(X^{(2)}\right) = \alpha + c\left(\tilde{F_2}\right) \ , \ \dots \ , \ \lambda_2 \left(X^{(p)}\right)  = \alpha + c\left(\tilde{F_p}\right).$$
So, $c\left(\tilde{F_1}\right) = c\left(\tilde{F_2}\right) = \dots = c\left(\tilde{F_p}\right).$
Also, for each $j$ in $[p]$ we have $\tilde{F_j} \subseteq X^{(j)}_0 \bigcup X^{(j)}_{{\rm index}\left(\tilde{F_j}\right)}$ and $\left|\tilde{F_j} \bigcap X^{(j)}_0\right| \leq s$.
Our aim is showing that
$\lambda_1 \left(X^{(1)}\right) , \lambda_1 \left(X^{(2)}\right) , \dots , \lambda_1 \left(X^{(p)}\right)$ are not pairwise distinct.
If $\lambda_1 \left(X^{(1)}\right) , \lambda_1 \left(X^{(2)}\right) , \dots , \lambda_1 \left(X^{(p)}\right)$ are pairwise distinct, then for any two arbitrary elements $i$
and $j$ in $[p]$ with $i < j$ we have
$$\begin{array}{lcccccccr}
	{\rm index}\left(\tilde{F_i}\right) \neq {\rm index}\left(\tilde{F_j}\right) & & {\rm and} & & 
	X^{(i)}_{{\rm index}\left(\tilde{F_i}\right)} \subseteq X^{(j)}_{{\rm index}\left(\tilde{F_i}\right)} & & {\rm and} & & X^{(j)}_{0} \subseteq X^{(i)}_{0} ;
\end{array}$$
and therefore,
$$\begin{array}[5]{llcll}
\tilde{F_i} \cap \tilde{F_j} 	&  \subseteq   &  \Big[ X^{(i)}_{0} \bigcup X^{(i)}_{{\rm index}\left(\tilde{F_i}\right)} \Big]   &    \bigcap    &   \Big[ X^{(j)}_{0} \bigcup X^{(j)}_{{\rm index}\left(\tilde{F_j}\right)} \Big]     \\
     &     &    &     &     \\
	&   \subseteq   &   \Big[ X^{(i)}_{0} \bigcup X^{(j)}_{{\rm index}\left(\tilde{F_i}\right)} \Big] &   \bigcap     &  \Big[ X^{(i)}_{0} \bigcup X^{(j)}_{{\rm index}\left(\tilde{F_j}\right)} \Big]      \\
	     &     &    &     &     \\
	&   =   &   X^{(i)}_{0}    &     \bigcup     &  \underbrace{ \Big[X^{(j)}_{{\rm index}\left(\tilde{F_i}\right)} \bigcap X^{(j)}_{{\rm index}\left(\tilde{F_j}\right)} \Big]  }_{=\varnothing } .    \\
\end{array}$$
We conclude that $\tilde{F_i} \cap \tilde{F_j}   \subseteq  \tilde{F_i} \cap   X^{(i)}_{0}   $;
and therefore, $\left| \tilde{F_i} \cap \tilde{F_j} \right|   \leq s.   $

\noindent This shows that $\left\{\tilde{F_1} , \tilde{F_2} , \dots , \tilde{F_p} \right\}$ is a hyperedge of
${\rm KG}^p({\mathcal F} , s)$ with $c\left(\tilde{F_1}\right) = c\left(\tilde{F_2}\right) = \dots = c\left(\tilde{F_p}\right)$;
which contradicts the properness of $c:V({\rm KG}^p({\mathcal F} , s)) \longrightarrow [t]$.
Accordingly, $\lambda_1 \left(X^{(1)}\right) , \lambda_1 \left(X^{(2)}\right) , \dots , \lambda_1 \left(X^{(p)}\right)$ are not pairwise distinct.

Now, since all the hypotheses of ${\mathbb Z_p}$-Tucker Lemma are satisfied, we obtain
$$\alpha + (m-\alpha )(p-1) \geq n;$$
and therefore,
$$\chi({\rm KG}^p({\mathcal F} , s))\ge \left\lceil \frac{{\rm ccd}^p({\mathcal F} , s )}{p-1}\right\rceil,$$
as desired.	
\end{proof}

\begin{remark}
	In the Proof of Theorem \ref{Saeed}, for defining the map
		$$\lambda: ({\mathbb Z_p} \cup \{0\}) ^ n \setminus \{(0,0, \dots , 0)\} \longrightarrow \mathbb Z_p\times [m]$$
	    $$\lambda(X)=(\lambda_1(X),\lambda_2(X))$$
	one could also regard an arbitrary $X=(X_1, X_2, \dots , X_p)$ in $({\mathbb Z_p} \cup \{0\}) ^ n \setminus \{(0,0, \dots , 0)\}$, and consider the following two separate cases:
	
	\noindent
	{\bf Case 1.} If there exist a hyperedge 
	$F$ in $E({\mathcal F})$
	and a number $i$ in $[p]$ for which $F \subseteq X_0 \cup X_i$ and $\left| F \cap X_0 \right| \leq s$, then $\left| F \setminus X_0 \right| = |F| -  \left| F \cap X_0 \right| \geq |F| - s >0$. Now, since $F\cap X_i \neq \varnothing$ while $F\cap X_j = \varnothing$ for all other numbers $j$ in $[p] \setminus \{i\}$, one observes that such a number $i$ would absolutely depend on $F$ and it is certainly unique with respect to $F$. So, it could be easily called ${\rm index} (F)$. Now, for determining $\lambda (X)$, let $\tilde{F}$ be the maximum of all such mentioned hyperedges $F$ with respect to the total ordering $\preccurlyeq $. Put $$\lambda (X) := \left(\omega ^{{\rm index}\left(\tilde{F}\right)} \ , \ \alpha + c \left( \tilde{F} \right) \right).$$
			
	\noindent		
	{\bf Case 2.} If The Case 1 does not happen, then on account of the notion of ${\rm ccd}^p ( {\mathcal F} , s)$, we have
	$$\left| X_0 \right| \geq {\rm ccd}^p ( {\mathcal F} , s);$$
	and also,
	$$\sum_{j=1}^p \left| X_j \right| = n - \left| X_0 \right| \leq n - {\rm ccd}^p ( {\mathcal F} , s) = \alpha.$$
	Let
	$X_i$ be the maximum of $\{X_1 , X_2 , \dots , X_p\}$ with respect to the total ordering $\preccurlyeq$. Now, put
    $$\lambda (X) := \left(\omega ^{i} \ , \ \sum_{j=1}^p \left| X_j \right|  \right).$$
\end{remark}

\def\cprime{$'$} \def\cprime{$'$}

\end{document}